\documentclass [12pt,a4paper,reqno]{amsart}
\usepackage{amsmath,amsthm,amscd}
\usepackage{amsfonts,amssymb}
\allowdisplaybreaks

\newcommand{\be}{\begin{equation}}
\newcommand{\ef}{\end{equation}}
\chardef\bslash=`\\ 

\newcommand{\G}{\Gamma}
\newcommand{\wt}{\widetilde}
\newcommand{\wh}{\widehat}
\newcommand{\ov}{\overline}
 \renewcommand{\sectionmark}[1]{}

\newcommand{\iy}{\infty}
\newcommand{\Be}{Beltrami}
\newcommand{\hol} {holomorphic}
\newcommand{\qc} {quasiconformal}

\newcommand{\ve}{\varepsilon}

\newcommand{\fc}{\frac}

 \usepackage{amsfonts}

\newcommand{\g}{\gamma}

\newcommand{\D}{\mathbb{D}}

\newcommand{\z}{\zeta}
\newcommand{\vp}{\varphi}
\newcommand{\hC}{\widehat{\mathbb{C}}}
\newcommand{\C}{\mathbb{C}}
\newcommand{\R}{\mathbb{R}}

\newcommand{\B}{\mathbf{B}}
\newcommand{\T}{\mathbf{T}}

\newcommand{\Belt} {\operatorname{Belt}}

\newcommand{\Teich} {\operatorname{Teich}}

\begin{document}

\title{On holomorphic contractibility of Teichm\"{u}ller spaces}

\author{Samuel L. Krushkal}

\begin{abstract} The problem of holomorphic contractibilty of the
Teichm\"{u}ller spaces $\T(0, n)$ of punctured spheres ($n > 4$)
arose in the 1970s in connection with solving algebraic
equations in Banach algebras. Recently it was solved by the author
in \cite{Kr2}.

In the present paper we improve the statement of Lemma 3 in \cite{Kr2}
and provide an alternate proof of holomorphic contractibility of
low dimensional Teichm\"{u}ller spaces.
\end{abstract}

\date{\today\hskip4mm ({\tt OnContract1.tex})}

\maketitle

\bigskip

{\small {\textbf {2010 Mathematics Subject Classification:} Primary:
30C55, 30F60; Secondary: 30F35, 46G20}

\medskip


\bigskip

\markboth{Samuel Krushkal}{On holomorphic contractibility of
Teichm\"{u}ller spaces} \pagestyle{headings}

\bigskip\bigskip\centerline
{\bf 1. PREAMBLE}

\bigskip\noindent{\bf 1.1}.
A complex Banach manifold $X$ is contractible to its point
$x_0$ if there exists a continuous map
$F: X \times [0, 1] \to X$ with $F(x, 0) = x$ and $F(x, 1) = x_0$
for all $x \in X$.
If the map $F$ can be chosen so that for every $t \in [0, 1]$ the map
$F_t: x \mapsto F(x, t)$ of $X$ to itself is holomorphic and $F_t(x_0) = x_0$,
then $X$ is called holomorphically contractible to $x_0$.

The problem of holomorphic contractibilty of
Teichm\"{u}ller spaces $\T(0, n)$ of the punctures spheres ($n > 4$)
arose in the 1970s in connection with solving the algebraic
equations in Banach algebras. It was caused by the fact that in the space
$\mathbb C^m, \ m > 1$, there are domains (even bounded), which are only
topologically but not holomorphically contractible (see \cite{Go}, \cite{Hi},
\cite{ZL1}, \cite{ZL2}).

Recently this problem was solved positively in \cite{Kr2}. There is established
that all space $\T(0, n), \ n > 4$, are holomorphically contractible.

The proof of Lemma 3 in this paper contains a wrongly assertion (which does not be used  there) that the map $s_m$ (giving the inclusion of the space $T(\Gamma_0)$ into  $T(\Gamma^m_0)$) is a section of the forgetful map $\chi_m: \ T(X^m_{a^0}) \rightarrow T(X_{a^0})$. Such sections do not exist if $n>6$.

In fact, there was only used in the proof that $s_m$ is an open holomorphic map (of a domain onto manifold), and the openness is preserved for the limit map $s = \lim_{m\to\infty} s_m$ which determines an $(n-3)$-dimensional complex submanifold $s(T(X_{a^0}))$ in the universal Teichm\"uller space $\T$.

\bigskip
In the present paper, we improve the statement of the indicated Lemma 3 (without
changing other arguments in \cite{Kr2} concerning the contractibility of
$\T(0, n)$). In the second part of the paper, we provide an alternate proof of holomorphic contractibility of low dimensional Teichm\"{u}ller spaces (of dimensions two and three), which has an independent interest in view of importance of the problem.
The underlying idea of this proof is different from \cite{Kr2}.

\bigskip\noindent{\bf 1.2}.
There are two Teichm\"{u}ller spaces of dimension two: the space $\T(0, 5)$
of the spheres with five punctures and the space $\T(1,2)$ of twice punctured
tori; these spaces are biholomorphically equivalent. Such spheres and tori are uniformized  by the corresponding Fuchsian groups $\Gamma$ and
$\Gamma^\prime$ so that $\Gamma$ is a subgroup of index two in $\Gamma^\prime$;
letting  $\T(0,5) = \T(\Gamma), \ \T(1, 2) = \T(\Gamma^\prime)$ ,
we have $\T(\Gamma^\prime) = \T(\Gamma)$.

In a similar way, the Teichm\"{u}ller spaces $\T(0, 6)$ of spheres with six punctures
and $\T(2, 0)$ of closed Riemann surfaces of genus $2$ also are biholomorphically equivalent, and in terms of the corresponding Fuchsian groups $\Gamma$ and
$\Gamma^\prime$  we have the same relation $\T(\Gamma^\prime) = \T(\Gamma)$.
We state:

\bigskip\noindent
{\bf Theorem 1}. {\em The spaces $\T(0, 5), \ \T(1, 2), \ \T(0, 6), \ \T(2, 0)$ are
holomorphically contractible.}

\bigskip
It remains the Teichm\"{u}ller space $\T(1, 3)$ of tori with three punctures also having dimension three; it will not be involved here.

\bigskip\bigskip
\centerline {\bf 2. UNDERLYING RESULTS}

\bigskip\noindent
{\bf 2.1}.
Consider the ordered $n$-tuples of points
\begin{equation}\label{1}
\mathbf a = (0, 1, a_1, \dots, a_{n-3}, \infty), \quad n > 4,
\end{equation}
with distinct $a_j \in \hC \setminus \{1, - 1, i\}$ and the corresponding
punctured spheres
$$
X_{\mathbf a} = \hC \setminus \{0, 1, a_1 \dots, a_{n-3}, \infty\},
\quad \hC = \C \cup \{\infty\},
$$
regarded as the Riemann surfaces of genus zero. Fix a collection
$\mathbf a^0 = (0, 1, a_1^0, \dots, a_{n-3}^0, 1, \infty)$ defining the base
point $X_{\mathbf a^0}$ of Teichm\"{u}ller space
$\T(0, n) = \T(X_{\mathbf a^0})$ of $n$-punctured spheres. Its points
are the equivalence classes $[\mu]$ of \Be \ coefficients from the
ball
$$
\Belt(\C)_1 = \{\mu \in L_\iy(\C): \ \|\mu\|_\iy < 1\},
$$
under the relation that $\mu_1 \sim \mu_2$ if the corresponding
\qc \ homeomorphisms
$w^{\mu_1}, w^{\mu_2}: \ X_{\mathbf a^0} \to X_{\mathbf a}$
(the solutions of the \Be \ equation
$\overline{\partial} w = \mu \partial w$ with $\mu = \mu_1, \mu_2$)
are homotopic on $X_{\mathbf a^0}$ (and hence coincide
in the points $0, 1, a_1^0, \dots, a_{n-3}^0, \iy$).
This models $\T(0, n)$ as the quotient space
$$
\T(0, n) = \Belt(\C)_1/\sim
$$
with complex Banach structure of dimension $n - 3$ inherited from
the ball $\Belt(\C)_1$.

Another canonical model of $\T(0, n) = \T(X_{\mathbf a^0})$ is
obtained using the uniformization of Riemann surfaces and the
holomorphic Bers embedding of Teichm\"{u}ller \ spaces.
Consider the upper and lower half-planes
$$
U  = \{z = x + i y: \ y > 0\}, \quad U^* = \{z \in \hC: \ y < 0\}
$$
The holomorphic universal covering map $h: \ U \to X_{\mathbf a^0}$
provides is a torsion free Fuchsian group $\G_0$ of the first kind
acting discontinuously on $U \cup U^*$, and the surface $X_{\mathbf a^0}$
is represented (up to conformal equivalence) as the quotient space
$U/\G_0$.
The functions $\mu \in L_\iy(X_{\mathbf a^0}) = L_\iy(\C)$ are lifted
to $U$ as the Beltrami $(-1, 1)$-measurable forms  $\wt \mu d\ov{z}/dz$
on $U$ with respect to $\G_0$ satisfying
$(\wt \mu \circ \g) \ov{\g^\prime}/\g^\prime = \wt \mu, \ \g \in
\G_0$ and form the corresponding Banach space $L_\iy(U, \G_0)$.
We extend these $\wt \mu$ by zero to $U^*$ and consider the unit
ball $\Belt(U, \G_0)_1$ of this space $L_\iy(U, \G_0)$. The corresponding
quasiconformal maps $w^{\wh{\mu}}$ are conformal  on the half-plane $U^*$,
and their {\it Schwarzian derivatives}
$$
S_w(z) = \Bigl(\frac{w^{\prime\prime\prime}(z)}{w^\prime(z)}\Bigr)^\prime
- \frac{1}{2} \Bigl(\fc{w^{\prime\prime}(z)}{w^\prime(z)}\Bigr)^2,
\quad w = w^{\wh{\mu}},
$$
fill a bounded domain in the complex $(n - 3)$-dimensional space $\B(\G_0) = \B(U^*, \G_0)$ of hyperbolically bounded holomorphic $\G_0$-automorphic forms of degree $- 4$
on $U^*$ (i.e.,  satisfy $(\vp \circ \g) (\g^\prime)^2 = \vp, \ \g \in \G_0$),
with norm
$$
\|\vp\|= \sup_{U^*} 4 y^2 |\vp(z)|.
$$
This domain models the {\it Teichm\"{u}ller space $\T(\G_0)$ of the group $\G_0$}.
It is canonically isomorphic to the space $\T(X_{\mathbf a^0})$ (and to  a
bounded domain in the complex Euclidean space $\C^{n-3}$).

The indicated map  $\wh{\mu} \to S_{w^{\wh{\mu}}}$ determines a holomorphic  a map $\phi_\T: \Belt(U, \G_0)_1 \to \B(\G_0)$; it has only local holomorphic sections.

Note also that $\T(\G_0) = \T \cap \B(\G_0)$, where $\T$ is the
universal Teichm\"{u}ller space (modelled as domain of the Schwarzian derivatives
of all univalent functions on $U^*$ admitting quasiconformal extension to $U$).

\bigskip\noindent
{\bf 2.2}. The collections (1) fill a domain $\mathcal D_n$ in
$\C^{n-3}$ obtained by deleting from this space the hyperplanes $\{z
= (z_1, \dots, z_{n-3}): \ z_j = z_l, \ j \ne l\}$, and with $z_1 =
0, z_2 = 1$.
This domain represents the Torelli space of the spheres
$X_{\mathbf a}$ and is covered by $\T(0, n)$, which is given by the following
lemma (cf. e.g., \cite{Ka}; \cite{Na}, Section 2.8).

\bigskip\noindent
{\bf Lemma 1}. {\em The holomorphic universal covering space of $\mathcal
D_n$ is the Teichm\"{u}ller \ space $\T(0, n)$. This means that for
each punctured sphere $X_{\mathbf a}$ there is a \hol \ universal
covering
$$
\pi_{\mathbf a}: \T(0, n) = \T(X_{\mathbf a}) \to \mathcal D_n.
$$
The covering map $\pi_{\mathbf a}$ is well defined by
$$
\pi_{\mathbf a} \circ \phi_{\mathbf a}(\mu) =
(0, 1, w^\mu(a_1), \dots, w^\mu(a_{n-3}), \infty),
$$
where $\phi_{\mathbf a}$ denotes the canonical projection of the ball
$\Belt(\C)_1$ onto the space $\T(X_{\mathbf a})$. }

\bigskip
Lemma 1 yields also that the truncated collections $\mathbf a_{*}
= (a_1, \dots, a_{n-3})$ provide the local complex coordinates on
the space $\T(0, n)$ and define its complex structure.

These coordinates are simply connected with the Bers local complex
coordinates on $\T(0, n)$ (related to basises of the tangent spaces
to $\T(0, n)$ at its points, see \cite{Be1}) via standard variation
of quasiconformal maps of $X_{\mathbf a} = U/\G_{\mathbf a}$
$$
\begin{aligned}
w^\mu(z) &= z - \fc{z(z-1)}{\pi}\iint\limits_\C
\fc{\mu(\zeta)}{\zeta(\zeta-1)(\zeta-z)} d\xi d\eta + O(\|\mu\|_\iy^2) \\
&= z - \fc{z(z - 1)}{\pi} \sum\limits_{\g \in \G_{\mathbf a}} \
\iint\limits_{U/\G_{\mathbf a}} \fc{\mu(\g
\zeta)|\g^\prime(\zeta)|^2}{\g \zeta(\g \zeta-1)(\g \zeta-z)} d\xi
d\eta + O(\|\mu\|_\iy^2).
\end{aligned}
$$
with uniform estimate of the ratio $O(\|\mu\|_\iy^2)/\|\mu\|_\iy^2$ on compacts
in $\C$ (see e.g., \cite{Kr1}).

\bigskip
It turns out that one can obtain the whole space $\T(X_{\mathbf a^0})$ using only the similar
equivalence classes $[\mu]$ of the Beltrami coefficients from the ball $\mu \in \Belt(U)_1$ (vanishing on $U^*)$, requiring that the corresponding
quasiconformal homeomorphisms $w^\mu$  are homotopic on the punctured sphere
$X_{\mathbf a^0}$. Surjectivity of this holomorphic quotient map
$$
\chi: \Belt(U)_1 \to \T(0, n),
$$
is a consequence of following interpolation result from \cite{CHMG}.

\bigskip\noindent
{\bf Lemma 2}. {\it Given two cyclically ordered collections of
points $(z_1, \dots, z_m)$ and $(\z_1, \dots, \z_m)$ on the unit circle
$S^1 = \{|z| = 1\}$, there exists a \hol \ univalent function $f$
in the closure of the unit disk $\D = \{|z| < 1\}$ such that
$|f(z)| < 1$ for $z \in \ov \D$ distinct from $z_1, \dots, z_m$,
and $f(z_k) = \z_k$ for all $k = 1, \dots, m$. Moreover, there exist
univalent polynomials $f$ with such an interpolation property. }

\bigskip
It follows that the function $f$ given by Lemma 2 actually is
holomorphic and univalent (hence, maps conformally) in a broader disk $\D_r, \ r > 1$.

First of all, $f^\prime(z) \ne 0$ on the unit circle $S^1$, Indeed,
were $f^\prime(z_0) = 0$ at some point $z_0 \in S^1$, then in its neighborhood
$f(z) = c_s (z - z_0)^s + O((z - z_0)^{s+1}) = c_s \zeta^s$, where $c_s \ne 0$
for some $s > 1$, which contradicts to injectivity of $f(z)$ on $S^1$. So $f$ is
univalent in some disk $\D_r = \{|z| < r\}, \ r > 1$.

Now, assuming, in the contrary, that $f$ is not globally univalent in any
admissible disk $\D_r$ with $r > 1$, one obtains the distinct sequences
$\{z_n^\prime\}, \ \{z_n^{\prime\prime}\} \subset \D_r$  with $f(z_n^\prime) = f(z_n^{\prime\prime})$ for any $n$,
whose limit points $z_0^\prime, \ z_0^{\prime\prime}$ lie on $S^1$.
Then also in the limit, $f(z_0^\prime) = f(z_0^{\prime\prime})$, which
in the case $z_0^\prime \ne z_0^{\prime\prime}$ contradicts to univalence
of $f$ on $S^1$ given by Lemma 2 and in the case
$z_0^\prime =  z_0^{\prime\prime} = z_0$ contradicts to local univalence of
$f$ in a neighborhood of $z_0$.

\bigskip
The interpolating function $f$ given by Lemma 2 is extended quasiconformally
to the whole sphere $\hC$ across any circle $\{|z| = r\}$ with $r > 1$ indicated
above. Hence, given a
cyclically ordered collection $(z_1, \dots, z_m)$ of points on
$S^1$, then for any ordered collection $(\z_1, \dots, \z_m)$ in
$\hC$, there exists a quasiconformal homeomorphism $\wh f$ of the sphere
$\hC$ carrying the points $z_j$ to $\z_j, \ j = 1, \dots, m$, and
such that its restriction to the closed disk $\ov \D$ is
biholomorphic on $\ov \D$.

Taking the quasicircles $L$ passing through the points
$\z_1, \dots, \z_m$ and such that each $\z_j$ belongs to an analytic subarc of $L$,
one obtains quasiconformal extensions of the interpolating function $f$, which are
biholomorphic on the union of $\ov{\D}$ and some neighborhoods of the initial points
$z_1, \dots, z_m \in S^1$. Now Lemma 1 provides quasiconformal extensions of $f$ lying in prescribed homotopy classes of homeomorphisms $X_{\mathbf z} \to
X_{\mathbf w}$.

\bigskip\noindent
{\bf 2.3}. Pick a space $\T(0, n) = \T(X_{\mathbf a^0})$ with $n \ge  5$
and let
$$
X_{\mathbf a^0}^\prime = X_{\mathbf a^0} \setminus \{a_{n-3}^0\} = U/\G_0^\prime.
$$
Due to the Bers isomorphism theorem, {\it the space $T(X_{\mathbf a^0}^\prime)$ is
biholomorphically isomorphic to the Bers fiber space }
$$
\mathcal F(0, n) := \mathcal F(T(X_{\mathbf a^0})
= \{(\phi_\T(\mu), z) \in \T(X_{\mathbf a^0}) \times \C: \ \mu \in
\Belt(U, \G_0^\prime)_1, \ z \in w^\mu(\D)\}
$$
{\it over} $\T(X_{\mathbf a^0})$  {\it with holomorphic projection} $\pi(\vp,
z) = \vp \ (\vp \in T(X_{\mathbf a^0})$ (see \cite{Be2}).

\bigskip
This fiber space is a bounded hyperbolic domain in $\B(\G_0) \times \C$ and represents
the collection of domains $D_\mu = w^\mu(U)$ (the universal covers of
the surfaces $X_{\mathbf a^0}$) as a holomorphic family  over the space $\T(0, n - 1) = T(X_{\mathbf a^0})$.

The indicated isomorphism between $\T(0, n + 1)$  and $\mathcal F(0, n)$ is
induced by the inclusion map \linebreak $j: \ \D_{*} \hookrightarrow
\D$ forgetting the puncture at $a_n^0$, via
 \be\label{2}
\mu \mapsto (S_{w^{\mu_1}}, w^{\mu_1}(\wh a_{n-3}^0)) \quad \text{with} \ \
\mu_1 = j_{*} \mu := (\mu \circ \wh j_0) \ov{\wh{j}^\prime}/\wh{j}^\prime,
\end{equation}
where $\wh{j}$ is the lift of $j$ to $U$ and $\wh a_{n-3}^0$ is one of the points
from the fiber over $a_n^0$ under the quotient map $U \to U/\Gamma_0$.

Note also that the holomorphic universal covering maps $h: \ U^* \to U^*/\Gamma_0$
and $h^\prime: \ U^* \to U^*/\Gamma_0^\prime$ (and similarly the corresponding covering maps in $U$) are related by
$$
j \circ h^\prime = h \circ \wh j,
$$
where $\wh j$ is the lift of $j$. This induces a
surjective homomorphism of the covering groups $\theta: \Gamma_0
\to \Gamma_0^\prime$ by
$$
 \wh j \circ \g = \theta(\g) \circ \g, \quad \g \in \Gamma_0^\prime
$$
and the norm preserving isomorphism $\wh j_{*}: \ \B(\Gamma_0)
\to \B(\Gamma_0^\prime)$ by
 \begin{equation}\label{3}
\wh j_{*} \vp = (\vp \circ \wh j) (\wh j^\prime)^2,
\end{equation}
which projects to the surfaces $X_{\mathbf a^0}$ and $X_{\mathbf
a^0}^\prime$ as the inclusion of the space $Q(X_{\mathbf a^0})$ of
holomorphic quadratic differentials corresponding to $\B(\Gamma_0)$ into the
space $Q(X_{\mathbf a^0}^\prime)$ (cf. \cite{EK}).

The Bers theorem is valid for Teichm\"{u}ller space $\T(X_0
\setminus \{x_0\})$ of any punctured hyperbolic Riemann surface
$X_0 \setminus \{x_0\}$ and implies that $\T(X_0 \setminus \{x_0\})$
is biholomorphically isomorphic to the Bers fiber space $ \mathcal
F(\T(X_0))$ over $\T(X_0)$.

\bigskip\noindent
{\bf 2.4}. The group $\Gamma_0^\prime$ uniformizing the surface
$X_{\mathbf a^0}$ acts discontinuously on the fiber space $\mathcal F(\Gamma_0)$ as a group of biholomorphic maps by
 \begin{equation}\label{4}
\g(\phi_\T(\mu), z) = (\phi_\T(\mu), \g^\mu z),
\end{equation}
where $\mu \in \Belt(U, \Gamma_0), \ z \in w^\mu(U)), \ \g \in \Gamma_0$,
and
$$
\g^\mu \circ w^\mu  = w^\mu \circ \g
$$
(see \cite{Be2}). The quotient space
$$
\mathcal V(0, n) := \mathcal V(\Gamma_0) = \T(0, n + 1)/\Gamma_0
$$
is called the {\it $n$-punctured Teichm\"{u}ller curve} and depends only on the analytic type
of the group $\Gamma_0$. The projection $\pi: \ \mathcal F(0, n) \to \T(0, n)$ induces a holomorphic projection
 \begin{equation}\label{5}
\pi_1: \ \ \mathcal V(0, n) \to \T(0, n).
\end{equation}
This curve also is a complex manifold with fibers $\pi^{-1}(x) = X_{\mathbf a}$.

Due to the deep Hubbard-Earle-Kra theorem \cite{EK}, \cite{Hu}, {\it the projections
$\mathcal V(0, n) \to \T(0, n)$ and (4) have no holomorphic sections for any $n \ge 7$} (more generally, for all spaces $\T(\Gamma)$
corresponding to groups $\Gamma$ without elliptic elements, excluding the spaces
$\T(1,2) \simeq \T(0, 5)$ and $\T(2, 0) \simeq \T(0, 6)$).
Such sections exist for groups $\Gamma$ containing elliptic elements.

In the exceptional cases of $\T(1,2)$ and $\T(2, 0)$, there is a group $\Gamma^\prime$ which contains $\Gamma $ as a subgroup of
index two. Then $\T(\Gamma^\prime) = \T(\Gamma), \ \mathcal F(\Gamma^\prime) =
\mathcal F(\Gamma)$, and the elliptic elements $\g \in \Gamma^\prime$ produce the indicated holomorphic sections $s$ as the maps
 \begin{equation}\label{6}
\phi_\T(\mu) \mapsto (\phi_\T(\mu), w^\mu(z_0)),
\end{equation}
where $z_0$ is a fixed point of $\g$ in the half-plane $U$.
These sections are called the {\it Weierstrass sections} (in view of their connection
with the Weierstrass points of hyperelliptic surface $U/\Gamma$).
We describe these sections following \cite{EK}.

Consider also the {\it punctured fiber space} $\mathcal F_0(\Gamma)$ to be the largest open dense subset of $\mathcal F(\Gamma)$ on which $\Gamma$ acts freely, and let
$$
\mathcal V^\prime(\Gamma) = \mathcal F_0(\Gamma)/\Gamma.
$$
For $\Gamma$ with no elliptic elements, the universal covering space for of
$\mathcal V^\prime(g, n) = \mathcal V^\prime(\Gamma)$ is $\T(g, n + 1)$.

If $\Gamma$ contains elliptic elements $\g$, then any  holomorphic section
$\T(\Gamma) \to \mathcal F(\Gamma)$ is determined by the maps (6) so that
$w^\mu(z_0)$ is exactly one fixed point of the corresponding map (4) in the
fiber $w^\mu(U)$. These  holomorphic sections take their values in the set
$\mathcal V(\Gamma) \setminus \mathcal V^\prime(\Gamma)$ and do not provide in general
sections of projection (5).

In the case of spaces $\T(1, 2)$ and $\T(2, 0)$, either of the corresponding
curves $\mathcal V(1, 2)$ and $\mathcal V(2, 0)$ has a unique biholomorphic self-map
$\g$ of order two which maps each fiber into itself. The fixed-point locus of $\g$ is a finite set of connected closed complex submanifolds of $\mathcal V^\prime(g, n)$,
and the restriction of the map (5) to one of these submanifolds is holomorphic map
onto $\T(0, n)$; its inverse is a Weierstrass section. The restriction of $\g$ to
each fiber is a conformal involution of the corresponding hyperelliptic Riemann
surface interchanging its sheets, and the fixed points of $\g$ are the Weierstrass
points on this surface.

\bigskip
In dimension one, there are three biholomorphically isomorphic Teichm\"{u}ller spaces $\T(1, 0), \ \T(1, 1)$ and $\T(0, 4)$ (see, e.g., \cite{Be2}, \cite{Pa}).
We shall use the last two spaces. Their fiber space $\mathcal F(0, 4))$ is isomorphic to $\T(0, 5)$.

As a special case of the Hubbard-Earle-Kra theorem \cite{EK}, \cite{Hu}, we have:

\bigskip\noindent
{\bf Lemma 3}. {\it
(a) The curve $\mathcal V(0, 4)$ has for any its point $x$ a unique holomorphic section $s$ with $s(\pi_1(x)) = x$.

(b) If $\dim \mathcal V(g, n)^\prime > 1$, the only curves $\mathcal V(1, 2)^\prime$ and $\mathcal V(2, 0)^\prime$ have holomorphic sections, which are their Weierstrass sections.  }

\bigskip
In particular, the curve $\mathcal V(2, 0)$ has six disjoint holomorphic sections corresponding to the Weierstrass points of hyperelliptic surfaces of genus two.

\bigskip\bigskip
\centerline{\bf 3. HOLOMORPHIC MAPS OF $\T(0, n)$ INTO UNIVERSAL TEICHM\"{U}LLER}
\centerline{\bf SPACE}

\bigskip\noindent
{\bf 3.1}.  The universal Teichm\"{u}ller space $\T = \Teich (U)$
is the space of quasisymmetric homeomorphisms of the unit circle
factorized by M\"{o}bius maps;  all Teichm\"{u}ller
spaces have their isometric copies in $\T$.

The canonical complex Banach structure on $\T$ is defined by
factorization of the ball of the Beltrami coefficients
$$
\Belt(U)_1 = \{\mu \in L_\iy(\C): \ \mu|U^* = 0, \ \|\mu\|_\iy < 1\}
$$
(i.e., supported in the upper-half plane),
letting $\mu_1, \mu_2 \in \Belt(U)_1$ be equivalent if the
corresponding \qc \ maps $w^{\mu_1}, w^{\mu_2}$ coincide on
$\ov{\mathbb R} = \mathbb R \cup \{\iy\} = \partial U^*$
(hence, on $\ov{U^*}$). Such $\mu$ and the corresponding maps
$w^\mu$ are called $\T$-{\it equivalent}. The equivalence classes
$[w^\mu]_\T$ are in one-to-one correspondence with the Schwarzian
derivatives $S_w$ in $U^*$, which fill a bounded domain in the space
$\B = \B(U^*)$ (see \textbf{2.1}).

The map  $\phi_\T: \ \mu \to S_{w^\mu}$ is holomorphic and descends to a
biholomorphic map of the space $\T$ onto this domain, which we will
identify with $\T$. As was mentioned above, it contains as complex
submanifolds the Teichm\"{u}ller spaces of all hyperbolic
Riemann surfaces and of Fuchsian groups.

We also define on this ball another equivalence relation, letting
$\mu, \ \nu \in \Belt(U)_1$ be equivalent if $w^\mu (a_j^0)
= w^\nu (a_j^0)$ for all $j$ and the homeomorphisms $w^\mu, \ w^\nu$
are homotopic on the punctured sphere $X_{\mathbf a^0}$. Let us call
such $\mu$ and $\nu$ {\it strongly $n$-equivalent}.

\bigskip\noindent
{\bf Lemma 4}. {\em If the coefficients $\mu, \nu \in \Belt(U)_1$
are $\T$-equivalent, then they are also strongly $n$-equivalent. }

\bigskip
The proof of this lemma is given in \cite{Ga}.

In view of Lemmas 1 and 4, the above factorizations of the ball
$\Belt(U)_1$ generate (by descending to the equivalence classes)
a holomorphic map $\chi$ of the underlying space $\T$ into
$\T(0, n) = \T(X_{\mathbf{a^0}})$.

This map is a split immersion (has local holomorphic sections),
which is a consequence, for example, of the following existence theorem from
\cite{Kr1}, which we present here as

\bigskip\noindent
{\bf Lemma 5}. {\it Let $D$ be a finitely connected domain on the Riemann
sphere $\hC$. Assume that there are a set $E$ of positive two-dimensional Lebesgue measure and a finite number of points
 $z_1, ..., z_m$ distinguished in $D$. Let
$\alpha_1, ..., \alpha_m$ be non-negative integers assigned to $z_1, ..., z_m$, respectively, so that $\alpha_j = 0$ if $z_j \in E$.

Then, for a sufficiently small $\ve_0 > 0$ and $\ve \in (0,
\ve_0)$, and for any given collection of numbers $w_{s j}, s
= 0, 1, ..., \alpha_j, \ j = 1,2, ..., m$, which satisfy the conditions
$w_{0 j} \in D$, \
$$
|w_{0j} - z_j| \le \ve, \ \ |w_{1 j} - 1| \le \ve, \ \ |w_{s j}| \le
\ve \ (s = 0, 1, \dots   a_j, \ j = 1, ..., m),
$$
there exists a quasiconformal automorphism $h$ of domain $D$ which is
conformal on $D \setminus E$ and satisfies
$$
h^{(s)}(z_j) = w_{sj} \quad \text{for all} \ s =0, 1, ..., \alpha_j,
\ j = 1, ..., m.
$$
Moreover, the Beltrami coefficient $\mu_h$ of $h$ on $E$ satisfies
$\|\mu_h\|_\iy \le M \ve$. The constants $\ve_0$ and $M$ depend only upon
the sets $D, E$ and the vectors $(z_1, ..., z_m)$ and $(\alpha_1, ..., \alpha_m)$.}

\bigskip\noindent{\bf 3.2}.
In fact, we have more, given by the following theorem which corrects the assertion
of Lemma 3 in \cite{Kr2}, as was mentioned in Preamble.

\bigskip\noindent
{\bf Theorem 2}. {\em The map $\chi$ is surjective and generates an open holomorphic
map $s$ of the space $\T(0, n) = \T(X_{\mathbf{a^0}})$ into the universal Teichm\"{u}ller
space $\T$ embedding $\T(0, n)$ into $\T$ as a $(n-3)$-dimensional submanifold. }

\bigskip\noindent
{\bf Proof}. The surjectivity of $\chi$ is a consequence of Lemma 2. To construct $s$,
take a dense subset
$$
e = \{x_1, \ x_2, \ \dots\} \subset X_{\mathbf a^0} \cap \R
$$
accumulating to all points of $\R$ and consider the punctured spheres
$X_{\mathbf a^0}^m = X_{\mathbf a^0} \setminus \{x_1, \dots, x_m\}$ with $m > 1$.
The equivalence relations on $\Belt(\C)_1$ for $X_{\mathbf a^0}^m$ and
$X_{\mathbf a^0}$ generate the corresponding
holomorphic map $\chi_m: \T(X_{\mathbf a^0}^m) \to \T(X_{\mathbf a^0})$.

Uniformizing the surfaces $X_{\mathbf a^0}$ and $X_{\mathbf a^0}^m$ by the
corresponding torsion free Fuchsian groups $\Gamma_0$ and $\Gamma_0^m$ of the first
kind acting discontinuously on $U \cup U^*$ and applying
to $U^*/\Gamma_0$ and $U^*/\Gamma_0^m$ the construction from Section \textbf{2.3}
(forgetting the additional punctures), one obtains similar to (3) the norm preserving isomorphism $\wh j_{m,{*}}: \ \B(\Gamma_0) \to \B(\Gamma_0^m)$ by
$$
\wh j_{m,{*}} \vp = (\vp \circ \wh j) (\wh j^\prime)^2,
$$
which projects to the surfaces $X_{\mathbf a^0}$ and $X_{\mathbf
a^0}^m$ as the inclusion of the space $Q(X_{\mathbf a^0})$ of
quadratic differentials corresponding to $\B(\Gamma_0)$ into the
space $Q(X_{\mathbf a^0}^m)$, and (since the projection $\eta_m$ has local holomorphic sections) geometrically this relation yields
a holomorphic embedding of the space $\T(\Gamma_0)$ into $\T(\Gamma_0^m)$ as
an $(n-3)$-dimensional submanifold. Denote this embedding by $s_m$.

To investigate the limit function for $m \to \infty$, we compose the maps $s_m$
with the canonical biholomorphic isomorphisms
$$
\eta_m: \ \T(X_{\mathbf a^0}^m) \to \T(\Gamma_0^m) = \T \cap
\B(\Gamma_0^m) \quad (m = 1, 2, \dots).
$$
Then the elements of $\T(\Gamma_0^m)$ are given by
$$
\wh s_m(X_{\mathbf a}) = \eta_m \circ s_m(X_{\mathbf a}),
$$
and this is a collection of the Schwarzians $S_{f^m}(z)$ corresponding
to the points $X_{\mathbf a}$ of $\T(X_{\mathbf a^0})$. So for any
surface $X_{\mathbf a}$, we have
 \begin{equation}\label{7}
\wh s_m(X_{\mathbf a}) = S_{f^m}(z).
\end{equation}
Each $\Gamma_0^m$ is the covering group of the universal cover $h_m:
\ U^* \to X_{\mathbf a_0^m}$, which can be normalized (conjugating
appropriately $\Gamma_0^m$) by $h_m(- i) = - i, \ h_m^\prime(- i) >
0$. Take its fundamental polygon $P_m$ obtained as the union of the
circular $m$-gon in $U^*$ centered at $z = -i$ with the zero angles
at the vertices and its reflection with respect to one of the
boundary arcs. These polygons increasingly exhaust the half-plane
$U^*$ from inside; hence, by the Carath\'{e}odory kernel theorem,
the maps $h_m$ converge to the identity map locally uniformly in
$U^*$.

Since the set of punctures $e$ is dense on $\R$, it completely
determines the equivalence classes $[w^\mu]$ and $S_{w^\mu}$  as
the points of the universal space $\T$. The limit function
$\wh s = \lim_{m\to \iy} \wh s_m$ maps $\T(X_{\mathbf a^0}) = \T(0, n)$
into the space $\T$ and also is canonically defined by the marked spheres
$X_{\mathbf a}$.

Similar to (7), the function $\wh s$ is represented as the Schwarzian of
some univalent function $f^n$ on $U^*$ with quasiconformal extension to $\hC$
determined by $X_{\mathbf a}$.
Then, by the well-known property of elements in the functional spaces with sup
norms, $\wh s$ is holomorphic also in $\B$-norm on $\T$.

Lemma 5 yields that $\wh s$ is a locally open map from $\T(X_{\mathbf a^0})$ to $\T$.
So the image $\wh s(\T(X_{\mathbf a^0}))$ is an $(n-3)$-dimensional complex
submanifold in $\T$ biholomorphically equivalent to $\T(\G_0)$.
The proof of Theorem 2 is completed.

\bigskip
The holomorphy property indicated above is based on the following lemma of Earle \cite{Ea}.

\bigskip\noindent
{\bf Lemma 6}. {\em Let $E, T$ be open subsets of complex Banach
spaces $X, Y$ and $B(E)$ be a Banach space of holomorphic functions
on $E$ with sup norm. If $\vp(x, t)$ is a bounded map $E \times T
\to B(E)$ such that $t \mapsto \vp(x, t)$ is holomorphic for each $x
\in E$, then the map $\vp$ is holomorphic.}

\bigskip Holomorphy of $\vp(x, t)$ in $t$ for fixed $x$ implies the
existence of complex directional derivatives
$$
\vp_t^\prime(x,t) = \lim\limits_{\z\to 0} \fc{\vp(x, t + \z v) -
\vp(x, t)}{\z} = \fc{1}{2 \pi i} \int\limits_{|\xi|=1} \fc{\vp(x, t
+ \xi v)}{\xi^2} d \xi,
$$
while the boundedness of $\vp$ in sup norm provides the uniform
estimate
$$
\|\vp(x, t + c \z v) - \vp(x, t) - \vp_t^\prime(x,t) c v\|_{B(E)}
\le M |c|^2,
$$
for sufficiently small $|c|$ and $\|v\|_Y$.

\bigskip\noindent{\bf 3.2}.
Now the desired holomorphic homotopy of $\T(0, n) = \T(X_{\mathbf a^0})$ into its
base point is constructed as follows.

Using the canonical embedding of $\T(0, n)$
in $\T$ via $\T(\G_0)$, we define on the space $\T(\G_0)$ a holomorphic homotopy
applying the maps
$$
W^\mu = \sigma^{-1} \circ w^\mu \circ \sigma, \quad \mu \in \Belt(U)_1;
\quad \sigma(\z) = i(1 + \z)/(1 - \z), \ \ \z \in \D,
$$
and $w_t^\mu(z) := w^\mu(z, t) = \sigma \circ W_t^\mu \circ \sigma^{-1}(z)$;
then
$$
S_{w^\mu}(\cdot, t) = t^2 S_{w^\mu}(\cdot)
=  t^{-2} (S_{W^\mu} \circ \ \sigma^{-1}) (\sigma^\prime)^{-2}.
$$
This point-wise equality determines holomorphic map $\eta(\vp, t) =
S_{w^\mu_t}: \ \T \times \D \to \T$  with
$\eta(\mathbf 0, t) = \mathbf 0, \ \eta(\vp, 0) = \mathbf 0, \ \eta(\vp, 1) = \vp$.
It is not compatible with the group $\G_0$; hence, there are images $\eta(\vp, t) = S_{w_t^\mu}$ which are located in $\T$ outside of $\T(\G_0)$.

Composition of $\eta(\vp, t)$ with maps $\chi$ and $s$ carries these images to the
points of the space $\T(0, n) = \T(X_{\mathbf a^0})$ and implies the desired
holomorphic homotopy of $\T(0, n)$.

\bigskip\bigskip
\centerline {\bf 4. HOLOMORPHIC CONTRACTIBILITY OF LOW DIMENSIONAL}
\centerline{\bf TEICHM\"{U}LLER SPACES (PROOF OF THEOREM 1)}

\bigskip\noindent
{\it $(a)$ Case $n = 5$ (dimension two)}. It is enough to establish holomorphic contractibility of the space $\T(0, 5) \simeq \mathbf F(0, 4)$ for the spheres
$$
X_{\mathbf a} = \hC \setminus \{0, 1, a_1, a_2, \iy\}.
$$
The fibers of $\T(0, 5)$ are the spheres with quadruples of punctures
$\{0, 1, a_1, \iy\}$.

We start with construction of the needed holomorphic homotopy of the space
$\T(0, 5)$ to its
base point $X_{\mathbf a^0}$ and first apply the assertion $(a)$ of Lemma 3 about holomorphic sections over $\T(0,4)$. It implies for any point
$$
x = (S_{w^{\mu_1}}, w^{\mu_1}(\wh a_{n-3}^0)) \in \T(0, 5)
$$
a unique holomorphic section $s: \T(0, 5) \to \T(0, 4)$ with $s(\pi_1(x)) = x$.
This section has a common point with each fiber $\pi^{-1}(x)$ over $\T(0, 4)$.

Since $\T(0,4)$ is (up to a biholomorphic equivalence) a simply connected
bounded Jordan domain $D \subset \C$ containing the origin, there is a holomorphic isotopy
$h(\zeta, t): D \times [0, 1] \to D$ with $h(\zeta, 0) = \zeta, \ h(z, 1) = 0$.
Using this isotopy, we define a homotopy $h_1(\vp, t)$ on $\T(0, 5)$, which carries
each point
$x = (S_{w^\mu}, w^\mu(\wh a_2^0)) \in \T(0, 5)$ to its image on the section $s$
passing from $x$, that is
  \begin{equation}\label{8}
h_1(\vp, w^\mu(\wh a_2^0)) = (h(\vp), \wt a_2), \quad \vp = S_{w^\mu}, \ \
\mu \in \Belt(\C)_1,
\end{equation}
where $\wt a_2$ is the common point of the fiber $h(\vp)$ and the selected section $s$.
The holomorphy of this homotopy in variables $x = (S_{w^\mu}, w^\mu(\wh a_2^0))$
for any fixed $t \in [0, 1]$ follows from Lemmas 1, 2 and the Bers isomorphism theorem.
The limit map
$$
h_1^*(x) = \lim\limits_{t\to 1} h_1(x, t),
$$
carries each fiber $w^\mu(U)$, to the initial half-plane $U$.

There is a canonical holomorphic isotopy
  \begin{equation}\label{9}
h_2(\zeta, t): \ U \times [0, 1] \to U
\end{equation}
of $U$ into its point corresponding to the origin of $\T(0, 5)$.
Now put $\mathbf h(x, t)$ to be equal
to $h_1(x, 2t)$ for $t \le 1/2$ and equal to $h_2(x, 2t - 1)$ for $x = \zeta \in U$
and $1/2 \le t \le 1$.

This function is holomorphic in $x \in \T(0, 5)$ for any fixed
$t \in [0, 1]$ and hence provides the desired holomorphic homotopy of the space $\T(0, 5)$ into its base point.

\bigskip\noindent
{\it $(b)$ Case $n = 6$ (dimension three)}. This case is more complicate.

We prescribe to each ordered sextuple $\mathbf a = \{0, 1, a_1, a_2, a_3, \iy\}$
of distinct points the corresponding punctured sphere
 \begin{equation}\label{10}
X_{\mathbf a} = \hC \setminus \{0, 1, a_1, a_2, a_3, \iy\}
\end{equation}
and the two-sheeted closed hyperelliptic surface $\wh X_{\mathbf a}$ of genus
two with the branch points $0, 1, a_1, a_2,
\newline
a_3, \iy$.
The corresponding Teichm\"{u}ller spaces $\T(0, 6)$ and $\T(2, 0)$ coincide
up to a natural biholomorphic isomorphism. Note also that the collections $\mathbf a = \{0, 1, a_1, a_2, a_3, \iy\}$ provide the local complex coordinates on
each from the spaces $\T(0, 6)$ and $\T(2, 0)$.

In view of the symmetry of hyperelliptic surfaces, it suffices to deal with
the Beltrami differentials $\mu d \ov{z}/dz$ on $\wh X_{\mathbf a}$, which are
compatible with
conformal involution $J_{\mathbf a}$ of $\wh X_{\mathbf a}$, hence, satisfy
$\mu(J_{\mathbf a} z) = \mu(z) J_{\mathbf a}^\prime/\ov{J_{\mathbf a}^\prime}$.
In other words, these $\mu$ are obtained by lifting to $\wh X_{\mathbf a}$ of the
Beltrami coefficients on $X_{\mathbf a}$. This extends Lemma 2 and its consequences
on holomorphy in the neighborhoods of the boundary interpolation points to
the corresponding two-sheeted disks on hyperelliptic surfaces.

\bigskip
Fix a base point of $\T(2, 0)$ determining a Fuchsian group $\Gamma$ for which $\T(\Gamma) = \T(2, 0)$.
The corresponding Teichm\"{u}ller curve $\mathcal V(2,0)$ is a $4$-dimensional complex analytic manifold with projection
$\pi_1: \ \mathcal V(2,0) \to \T(2, 0)$ onto $\T(2, 0)$ such that for every $\vp \in  \T(2,0)$ the fiber $\pi_1^{-1}(\vp)$ is a hyperelliptic surface, determined by $\vp$
(see \textbf{2.4}).

Due to assertion $(b)$ of Lemma 3, this curve has for any point
$$
\wh X_{\mathbf a} = (S_{w^{\mu_1}}, w^{\mu_1}(\wh a_{n-3}^0)) \in \T(2, 0))
$$
six distinct holomorphic sections $\wh s_1, \dots, \wh s_6$, corresponding to the Weierstrass points of the surface $X_{\mathbf a}$, with $\wh s_j(\pi_1(X_{\mathbf a})) = X_{\mathbf a}$,
and either from these sections has one common point with every fiber over $\T(2, 0)$.
We lift these sections to the Bers fiber space $\mathcal F(\Gamma)$ covering $\mathcal V(2,0)$.

As was mentioned in \textbf{2.4}, these sections are generated by the space $\mathcal F(\Gamma^\prime) = \mathcal F(\Gamma)$ corresponding to the extension $\Gamma^\prime$ of $\Gamma$, for which $\Gamma$ is a subgroup of index two.
Every section $\wh s_j$ acts on $\T(\Gamma^\prime)$ via (6), where $z_0$ is now the corresponding Weierstrass point of hyperelliptic surface $\wh X_{\mathbf a}$, and
$\wh s_j$ is compatible with the action (2) of the Bers isomorphism.

Thus each $\wh s_j$ descends to a holomorphic map $s_j: \ \T(0, 6) \to \mathcal V(0, 6)$ of the underlying space $\T(0, 6)$ for the punctured spheres (10).
We choose one from these maps and denote it by $s$.

The features of sections $\wh s_j$ provide that the descended map $s$ also determines  for each point $z_0 \in X_{\mathbf a}$ its unique image on every fiber $w^\mu(X_{\mathbf a})$ with $\mu \in \Belt(X_{\mathbf a})_1$, and this image is the point $w^\mu(z_0)$.

The next preliminary construction consists of embedding the space $\T(0, 5)$ into $\T(0, 6)$, using the forgetting map (3). Its image $j_{*} \T(0, 5)$ is a connected submanifold in $\T(0, 6)$, and the corresponding fibers of the curve $\mathcal V(0, 6)$ over the points $j_{*} \vp \in j_{*} \T(0, 5)$ are the surfaces $w^{j_{*} \mu}(X_{\mathbf a})$ with $j_{*} \mu(z) = \mu(\wh j(z)) \wh j^\prime(z)/\ov{\wh j^\prime(z)}$. The covering domains $w^{j_{*} \mu}(U)$ over these surfaces fill a submanifold $\wt \T(0, 7) \subset \T(0, 7)$, which is is biholomorphically equivalent to the space $\T(0, 6)$.

Using the biholomorphic equivalence of the space $\T(0, 5)$ to its image $j_{*} \T(0, 5)$
in $\T(0, 6)$, we carry over to $j_{*} \T(0, 5)$ the result of the previous step
$(a)$ on holomorphic  contractibility of $\T(0, 5)$, which provides a holomorphic homotopy
 \begin{equation}\label{11}
h(j_{*} \vp, t): \ j_{*} \T(0, 5) \times [0, 1] \to j_{*} \T(0, 5) \quad \text{with} \ \
h(j_{*} \vp, 0) = j_{*} \vp, \ \ h(j_{*} \vp, 1) = \mathbf 0
\end{equation}
(here $\mathbf 0$ stands for the origin of $j_{*} \T(0, 5)$).

Now we may construct the desired holomorphic homotopy of $\T(0, 6)$ contracting
this space to its origin.

First, regarding $\T(0, 6)$ as the Bers fiber space $\mathcal F(0, 5)$ over $\T(0, 5)$
(whose fibers are the covers of surfaces
$X_{\mathbf a^\prime}$ with collections of five punctures $\mathbf a^\prime = (0, 1, a_1, a_2, \iy)$),
we apply the homotopy (11) and define for any pair $x = (j_{*} \vp, z)$ with
$\vp \in \T(0, 5)$ and $z \in X_{\mathbf a}$ the map
 \begin{equation}\label{12}
\wt h_1((j_{*} \vp, z), t) = (h(j_{*} \vp), t), w_t^{j_{*} \mu} (z)),
\quad \vp \in \T(0, 5),
\end{equation}
noting that the image point $w_t^{j_{*} \mu} (z)$ is uniquely determined on the surface $w^{h(j_{*} \mu)}(X_{\mathbf a})$ by the map $s$, as was indicated above.

The pairs $(j_{*} \vp, z)$ are located in the space $\mathcal F(0, 6)$
and fill there a three-dimentional submanifold $\wt \T(0, 6)$ biholomorphically equivalent to $\T(0, 6)$.

The homotopy (12) is well defined on $\wt \T(0, 6) \times [0, 1]$ and contracts
the set $\wt \T(0, 6)$ into the fiber $\wt U$ over the base point.
It is holomorphic with respect to the space variable $x = (j_{*} \vp, z)$
for any fixed $t \in [0, 1]$ and continuous in both variables.

In view of biholomorphic equivalence of $\wt \T(0, 6)$ to $\T(0, 6)$, (12)
generates a holomorphic homotopy $h_1(x, t)$ of the space $\T(0, 6)$
onto the initial fiber (half-plane) $U$ over the origin of $\T(0, 5)$.

It remains to combine this homotopy $h_1$ with the additional homotopy (9) of $U$ into 
its point corresponding to the origin of $\T(0, 6)$.
This provides the desired homotopy $h$ and completes the proof of Theorem 1.

\bigskip\bigskip
{\small\em{
\leftline{Department of Mathematics, Bar-Ilan
University, 5290002 Ramat-Gan, Israel}
\leftline{and Department of Mathematics, University of Virginia, Charlottesville, VA 22904-4137, USA}}}

\end{document}